\documentclass[12pt]{article}

\usepackage[english]{babel}
\usepackage{amssymb,amsmath}
\usepackage{hyperref}

\newcommand{\footremember}[2]{%
    \footnote{#2}
    \newcounter{#1}
    \setcounter{#1}{\value{footnote}}%
}

\newtheorem{theorem}{Theorem }[section]

\newtheorem{corollary}[theorem]{Corollary}
\newtheorem{proposition}[theorem]{Proposition}

\newcommand{\qed}{\hspace*{\fill}$\Box$}

\def\min{\mathop{\mathrm{min}}\nolimits}

\title{ New necessary conditions for  (negative) Latin square type partial difference sets in abelian groups}

\author{Zeying Wang \footremember{MTU}{Michigan Technological University} \footnote{zeying@mtu.edu}}

\date{}

\begin{document}

\maketitle 

\abstract{Partial difference sets (for short, PDSs) with parameters ($n^2$, $r(n-\epsilon)$, $\epsilon n+r^2-3\epsilon r$, $r^2-\epsilon r$) are called {\it Latin square type} ({respectively {\it negative Latin square type}) PDSs if $\epsilon=1$ (respectively $\epsilon=-1$).  In this paper, we will give restrictions on the parameter $r$ of a (negative) Latin square type partial difference set in an abelian group of non-prime power order. As far as we know no previous general restrictions on $r$ were known. Our restrictions are particularly useful when $a$ is much larger than $b$.  As an application, we show that if there exists an abelian negative Latin square type PDS with parameter set $(9p^{4s},  r(3p^{2s}+1),-3p^{2s}+r^2+3r,r^2+r)$, $1 \le r \le \frac{3p^{2s}-1}{2}$, $p\equiv 1 \pmod 4$  a prime number and $s$ is an odd positive integer, then there are at most three possible values for $r$. For two of these three $r$ values, J. Polhill gave constructions in 2009 \cite{Polhill2009}.  }\\

\bigskip

\section{Introduction}

Let $G$ be a finite abelian group of order $v$, and let $D\subseteq G$ be a subset of size $k$. We say $D$ is a $(v,k,\lambda,\mu)$-{\it partial difference set} (PDS) in $G$ if the expressions $gh^{-1}$, $g$, $h\in D$, $g\neq h$, represent each non-identity element in $D$ exactly $\lambda$ times, and each non-identity element of $G$ not in $D$ exactly $\mu$ times. If we further assume that $D^{(-1)}=D$  (where $D^{(s)}=\{g^s:g\in D\}$) and  $e \notin D$ (where $e$ is the identity element of $G$), then $D$ is called a {\it regular} PDS.  A regular PDS is called {\it trivial} if $D\cup\{e\}$ or $G\setminus D$ is a subgroup of $G$. The condition that $D$ be regular is not a very restrictive one, as  $D^{(-1)}=D$ is automatically fulfilled whenever $\lambda\neq\mu$, and $D$ is a PDS if and only if $D\cup\{e\}$ is a PDS.

The {\it Cayley graph over $G$ with connection set $D$}, denoted by Cay($G$, $D$), is the graph with the elements of $G$ as vertices, and in which two vertices $g$ and $h$ are adjacent if and only if $gh^{-1}$ belongs to $D$.  When the connection set $D$ is a regular PDS, Cay($G$, $D$) is a strongly regular graph.  The importance of regular PDSs lies in the fact that they are equivalent to strongly regular graphs with a sharply transitive automorphism group. 

In  \cite{MA94b}, S.L. Ma describes several families of PDSs. Among these are PDSs with parameters ($n^2$, $r(n-\epsilon)$, $\epsilon n+r^2-3\epsilon r$, $r^2-\epsilon r$) for $\epsilon=\pm 1$. When $\epsilon=1$, the PDS is called a {\it Latin square type partial difference set}, and when $\epsilon=-1$, the PDS is called a {\it negative Latin square type partial difference set}. These types of PDSs have received a lot of attention, see for example \cite{Bailey}, \cite{Davis2004}, \cite{Davis2006}, \cite{Polhill2008}, \cite{Polhill2009_JCD}. 

Let $D$ be a regular $(v,k,\lambda, \mu)$-PDS, then $G\setminus D$ is a $(v, v-k, v-2k+\mu, v-2k+\lambda)$-PDS, and $(G\setminus D)\setminus\{e\}$ is a $(v, v-k-1, v-2k-2+\mu, v-2k+\lambda)$-PDS, and we call $(G\setminus D)\setminus\{e\}$ {\it the complement of $D$}. It is easy to verify that the complement of a (negative) Latin square PDS is again a (negative) Latin square type PDS (with the new parameter $r'=n-r+\epsilon$).  Thus when we study the existence of a (negative) Latin square type partial difference set (up to complements), we can add the restriction $1 \le r \le \frac{n+\epsilon}{2}$.  As far as we are aware, no other general restrictions on the parameter $r$ are known.

\medskip

In this article we will give restrictions on the parameter $r$ of Latin square type and negative Latin square type PDSs in abelian groups of order $a^2 b^2$, where $\gcd(a,b)=1$, $a>1$,  and $b$ is an odd positive integer $\ge 3$. Before stating our main result in the next section we provide an intuitive description of it here.  As mentioned before, $r$ could theoretically take any value in the interval $[1, \frac{ab+\epsilon}{2}]$. We divide this interval in subintervals of length approximately $a$. We then show that within each of these subintervals $r$ can only take values in two subintervals of length approximately $b$. When $a$ is less than approximately $2b$ this does not say anything, but when $a$ is much larger than $b$ this reduces the possible values $r$ can take significantly, as it reduces the number of possible values for $r$ to approximately $b^2$.

Throughout this paper, we will use the following standard notation: $\beta=\lambda-\mu$ and  $\Delta=\beta^2+4(k-\mu)$.

\section{Main result and its proof}

Below we cite a result on ``sub-partial difference sets"  that was discovered in \cite{MA94a}. We quote this result in the form given in \cite{MA94b} (Theorem 7.1). 

\begin{proposition}{\rm \cite{MA94a}}\label{Ma1}
Let D be a nontrivial regular $(v,\,k,\,\lambda,\,\mu)$-PDS  in an abelian group G. Suppose $\Delta$ is a  perfect square.  Let  $N$ be a subgroup of $G$ such that $\gcd(\left|N\right|, \left|G\right|/\left|N\right|)=1$ and $\left|G\right|/\left|N\right|$ is odd. Let 
$$\pi:=\gcd(|N|,\;\sqrt{\Delta}) \quad {\mbox and} \quad \theta:=\lfloor\frac{\beta+\pi}{2\pi}\rfloor.$$

Then  $D_1=N\,\cap D$ is a (not necessarily non-trivial)  regular $(v_1,\,k_1,\,\lambda_1,\,\mu_1)$-PDS in $N$ with
$$v_1=|N|, \; \beta_1=\lambda_1-\mu_1=\beta-2\theta \pi,  \; \Delta_1=\beta_1^2+4(k_1-\mu_1)=\pi^2,$$
and
 $$k_1=|N\cap D|=\frac{1}{2}\left[ |N| +\beta_1\pm \sqrt{(|N|+\beta_1)^2-(\Delta_1-\beta_1^2)(|N|-1)}  \right].$$

\end{proposition}

\medskip

{\bf Note:} In Proposition \ref{Ma1}, if $\Delta_1=|N|$, we have $$k_1=\frac{1}{2}\left[ |N| +\beta_1\pm (\beta_1+1)\sqrt{|N|}  \right].$$

Now we state our main theorem:

\begin{theorem}\label{Main}
Let $G$ be an abelian group of order $a^2b^2$, where $\gcd(a,b)=1$, $a>1$, and $b$ is an odd positive integer $\ge 3$. Let $D$ be an ($a^2 b^2$, $r(ab-\epsilon)$, $\epsilon ab+r^2-3\epsilon r$, $r^2-\epsilon r$)-PDS in $G$, where $$t \le \frac{\epsilon ab-2\epsilon r+a}{2a}<t+1$$
for some integer $t$. Then we have the following results:
\begin{itemize}
\item[{\rm (i)}]when $\epsilon=1$, that  is, when $D$ is a Latin square type PDS, we have 
either $$\frac{(a+1)b(b+1-2t)}{2(b+1)}-(b-1)\le r \le \frac{(a+1)b(b+1-2t)}{2(b+1)}$$
or
$$\frac{( a-1) b ( b-1 - 2 t)}{2(b-1)}\le r \le \frac{(a-1) b ( b-1 - 2 t)}{2(b-1)}+(b+1).$$

\item[{\rm (ii)}] when $\epsilon=-1$, that is, when $D$ is a negative Latin square type PDS, we have
either 
$$\frac{(a+1)b(b-1+2t)}{2(b-1)}-(b+1)\le r \le \frac{(a+1)b(b-1+2t)}{2(b-1)}$$
or
$$\frac{(a-1)b(b+1+2t)}{2(1+b)}\le r \le \frac{(a-1)b(1+b+2t)}{2(1+b)}+(b-1).$$
\end{itemize}
\end{theorem}

Note that for fixed $t$ the condition $t \le \frac{\epsilon ab-2\epsilon r+a}{2a}<t+1$ bounds $r$ to an interval of approximately length $a$, whereas the conclusion of the theorem bounds $r$ to two subintervals of approximately length $b$.

The proof of our main theorem is an application of Proposition \ref{Ma1} combined with a variance argument.

\medskip

{\bf Proof of Theorem \ref{Main}:} As $D$ is an ($a^2 b^2$, $r(ab-\epsilon)$, $\epsilon ab+r^2-3\epsilon r$, $r^2-\epsilon r$)-PDS in $G$,  we have  $\Delta=(\lambda-\mu)^2+4(k-\mu)=(\epsilon a b-2\epsilon r)^2+4(rab-r^2)=a^2b^2$.   Let $N$ be a subgroup of $G$ of order $a^2$. By Proposition \ref{Ma1}, we have $\pi=\gcd(|N|, \sqrt{\Delta})=a$, $\Delta_1=\pi^2=a^2=|N|$. Furthermore, $\beta_1=\beta-2\theta \pi$, and  $\theta=\lfloor \frac{\beta+\pi}{2\pi}\rfloor$. When 
$$t \le \frac{\beta+\pi}{2\pi}=\frac{\epsilon ab-2 \epsilon r+\pi}{2 \pi}<t+1,$$
we have $\theta=t$, and $\beta_1=\beta-2t\pi=\epsilon ab -2\epsilon r-2ta$.

Since $\Delta_1=|N|$,  by Proposition \ref{Ma1}, we have
\begin{eqnarray}
k_1&=& |N \cap D|=\frac{1}{2}\big(\,|N|+\beta_1\pm (\beta_1+1)\sqrt{|N|}\,\big)\nonumber\\
        &=& \frac{1}{2}\big((a^2+\epsilon ab -2\epsilon r -2t a)\pm(\epsilon ab-2\epsilon r -2ta+1)a\big). \label{size_k1} \end{eqnarray}

Assume that $G= N \times H$, where $H$ is a subgroup of $G$ of order $b^2$. Let $h_1$, $h_2$, $\cdots$, $h_{b^2-1}$ be the non-identity elements of $H$, let $\mathcal{B}_{h_i}=Nh_i\cap D$, and $B_i=|\mathcal{B}_{h_i}|$, $i=1, 2, \cdots, b^2-1$.

 One easily sees that $xy^{-1}$ ($x$, $y\in D$, $x \neq y$) belongs to $N$  if and only if either both $x$ and $y$ belong to $N$ or both belong to $\mathcal{B}_{h}$ for some non-identity element $h\in H$. By double-counting, we have
\begin{equation}\label{basic}
\sum_{i=1}^{b^2-1} B_i(B_i-1)+k_1(k_1-1)=k_1 \lambda+(a^2-1-k_1)\mu.
\end{equation}
It is easy to see that $\sum B_i=k-k_1$. Combining this information with Equation (\ref{basic}), we obtain the following system of equations:
\begin{equation}\label{main_formula}
\left\{\begin{array}{rl}
\sum B_i&=k-k_1,\\
\sum B_i^2&=k-k_1^2+k_1(\lambda-\mu)+(a^2-1)\mu.\end{array}\right.
\end{equation}

\begin{itemize}
\item[Case 1:] When $D$ is a Latin square type PDS, that is, when $\epsilon=1$, by (\ref{size_k1})  we have 
$k_1=\frac{1}{2}\left((a^2+ab -2 r -2t a)\pm(ab-2 r -2ta+1)a\right)$.
\begin{itemize}
\item[(i)] If $k_1=\frac{1}{2}\left(a^2+ab-2r-2ta)+(ab-2r-2ta+1)a\right)$, by (\ref{main_formula})
 the variance equals
\begin{eqnarray*}
&&(b^2-1)\sum B_i^2-(\sum B_i)^2\\
&&=-\frac{1}{4}a^2 (2+b+ab-b^2+ab^2-2 r-2 b r-2 bt-2 ab t)\\
&& *(b+ab+b^2+ab^2-2 r-2 b r-2 b t -2 ab t), \end{eqnarray*}
and hence is non-negative if and only if
$$\frac{(a+1)b(b+1-2t)}{2(b+1)}-(b-1)\le r \le \frac{(a+1)b(b+1-2t)}{2(b+1)}.$$

\item[(ii)] If  $k_1= \frac{1}{2}[(a^2+ab-2r-2t a)-(ab-2r-2t a+1)a]$, by (\ref{main_formula}) 
the variance equals
\begin{eqnarray*}
&&(b^2-1)\sum B_i^2-(\sum B_i)^2\\
&&=-\frac{1}{4} a^2(b - ab - b^2 + ab^2 + 2 r - 2 b r + 
    2 b t - 2 ab t) \\
&&*(-2 + b - ab + b^2 + 
    ab^2 + 2 r - 2 b r + 2 b t - 2 ab t), 
\end{eqnarray*}

and hence is non-negative if and only if
$$\frac{( a-1) b ( b-1 - 2 t)}{2(b-1)}\le r \le \frac{(a-1) b ( b-1 - 2 t)}{2(b-1)}+(b+1)$$
\end{itemize}

\end{itemize}

\item[Case 2:] When $D$ is a  negative Latin square type PDS, that is, when $\epsilon=-1$, by (\ref{size_k1}) we have 
$k_1=\frac{1}{2}\left((a^2-ab+2 r -2t a)\pm(-ab+2r -2ta+1)a\right)$.
 \begin{itemize}
\item[(i)] If $k_1=\frac{1}{2}[(a^2-ab+2r-2ta)+(-ab+2r-2ta+1)a]$, by (\ref{main_formula})
 the variance equals
\begin{eqnarray*}
&&(b^2-1)\sum B_i^2-(\sum B_i)^2\\
&&=-\frac{1}{4}a^2 (2-b-ab-b^2+ab^2+2r-2b r+2bt+2ab t)\\
&&*(-b-ab+b^2+ab^2+2r-2br+2bt+2abt), 
\end{eqnarray*}
and hence is non-negative if and only if
$$\frac{b(1+a)(b-1+2t)}{2(b-1)}-(b+1)\le r \le \frac{b(1+a)(b-1+2t)}{2(b-1)}.$$

\item[(ii)] If $k_1=\frac{1}{2}[(a^2-ab+2r-2ta)-(-ab+2r-2ta+1)a]$, by (\ref{main_formula}) the variance equals
\begin{eqnarray*}
&&(b^2-1)\sum B_i^2-(\sum B_i)^2\\
&&=-\frac{1}{4}a^2(-b+ab-b^2+ab^2-2r-2b r-2b t+2abt)\\
&&*(-2-b+ab+b^2+ab^2-2r-2br-2bt+2abt),
\end{eqnarray*}
and hence is non-negative if and only if
$$\frac{b(-1+a)(1+b+2t)}{2(1+b)}\le r \le \frac{b(-1+a)(1+b+2t)}{2(1+b)}+(b-1).$$
\end{itemize}

\qed

\section{Application  to Latin square type PDSs}

Letting $\epsilon=1$  and $t=(b-1)/2$  in Theorem \ref{Main}, we obtain the following Corollary:

\begin{corollary}\label{Latin_cor}
Let $G$ be an abelian group of order $a^2b^2$, $\gcd(a,b)=1$, $a>1$, and $b$ an odd positive integer $\ge 3$. Let $D$ be an ($a^2 b^2$, $r(ab-1)$, $ab+r^2-3r$, $r^2-r$)-Latin square type PDS in $G$ with $1 \le r \le a $.  Then we have
either $$1\le r \le  b+1,$$
or
$$\frac{(a+1)b}{(b+1)}-(b-1)\le r \le \frac{(a+1)b}{(b+1)}.$$
\end{corollary}

\medskip
When $b$ is the smallest prime power of $ab$, there are constructions satisfying the restriction $1 \le r \le b+1$ of Corollary \ref{Latin_cor}, see \cite{MA94b}. We cite the result below:

\begin{theorem} {\rm (\cite{MA94b})}
Let $n=p_1^{a_1} p_2^{a_2} \cdots p_s^{a_s}$ where $p_1$, $p_2$, $\cdots$, $p_s$ are distinct primes. Then there exists an abelian regular $(n^2, r(n-1), n+r^2-3r, r^2-r)$-PDS whenever $1 \le r \le \min\{p_i^{a_i}+1\}$.
\end{theorem}

\section{Applications to negative Latin square type PDSs}

In this section, we will apply our main result to negative Latin square type PDSs. In particular, we will show that in an abelian group $G$ of order $9m^4$,  where $m$ is an odd positive integer and $\gcd(3, m)=1$,  if there exists a negative Latin square type PDS with the parameter set $(9m^4,  r(3m^2+1),-3m^2+r^2+3r,r^2+r)$ and $1 \le r \le \frac{3m^2-1}{2}$, then there are at most 6 possible values for $r$. Furthermore, if $m=p^s$, where $p\equiv 1 \pmod 4$ is a prime number and $s$ is an odd positive integer, then there are at most 3 possible values for $r$. For two of these three $r$ values, Polhill gave constructions in 2009 \cite{Polhill2009}. Note that without our Main Theorem \ref{Main},  the possible number of values $r$ could hypothetically take in such groups grows unboundedly as $m$ or $p$ grows.

\medskip
Let $G$ be an abelian group of order $9m^4$, where $m$ is an odd positive integer and $\gcd(3,m)=1$.  Let $D$ be a $(9m^4, r(3m^2+1),-3m^2+r^2+3r,r^2+r)$ negative Latin square type PDS in $G$. Let $\epsilon=-1$, $t=-1$, $a=m^2$, and $b=3$ in Theorem \ref{Main}. It follows that if $0 \le r <m^2$, we have $\frac{3m^2-3}{4} \le r \le \frac{3m^2+5}{4}$. Let $\epsilon=-1$, $t=0$, $a=m^2$, and $b=3$ in Theorem \ref{Main}.  It follows that if $m^2\le r <2m^2$, then either $\frac{3m^2-5}{2}\le r \le \frac{3m^2+3}{2}$ or $\frac{3m^2-3}{2}\le r \le \frac{3m^2+1}{2}$.
\medskip

In what follows we will need the following result of Leung, Ma, and Schmidt \cite{Leung2008}.

\begin{theorem}\rm {(\cite{Leung2008})}\label{Ma2}
Let $p$ be a prime and let $D$ be a regular PDS in an abelian group $G=H \times P$ where $P$ is a $p$-group and $|H|$ is not divisible by $p$. Suppose that the parameter $\Delta$ is a square. Let $h \in H$. Then
$$
|D\cap Ph|\equiv\left\{ \begin{array}{cc}
1 & \pmod{p-1}\; \mbox{if $h \in D$}\\
0 & \pmod{p-1} \;\mbox{if $h \notin D$}
\end{array}\right.$$
\end{theorem}

\medskip

\begin{theorem}\label{NLS_1}
Let $G$ be an abelian group of order $9p^{4s}$, where $p$ is a prime number  $\ge 5$, and $s$ is a positive integer.  Then there does not exist a  $(9p^{4s},  r(3p^{2s}+1),-3p^{2s}+r^2+3r,r^2+r)$-PDS in $G$ with $r=(3p^{2s}-5)/2$.
\end{theorem}

{\bf Proof:}  Assume by way of contradiction that $D$ is a $(9p^{4s},  r(3p^{2s}+1),-3p^{2s}+r^2+3r,r^2+r)$-PDS in $G$ with $r=(3p^{2s}-5)/2$.  Let $N$ be the subgroup of $G$ of order $p^{4s}$. Applying Proposition \ref{Ma1} to the subgroup $N$,  we obtain $\pi=p^{2s}$, $\theta=0$, $\beta_1=-5$, and 
$$k_1= |N \cap D|=\frac{p^{4s}-4p^{2s}-5}{2}\;\; \mbox{or} \;\;\frac{p^{4s}+4p^{2s}-5}{2}.$$

Let $G=N\times H$, where $H$ is the subgroup of $G$ of order 9. As in the proof of Theorem \ref{Main}, let $h_1$, $h_2$, $\cdots$, $h_8$ be all non-identity elements of $H$, and let $\mathcal{B}_{h_i}=Nh_i \cap D$,  $B_i=|\mathcal{B}_{h_i}|$, $i=1, 2, \cdots, 8$. Now we discuss the existence of $D$ in two cases based on the possible size of $N \cap D$. 

\begin{itemize}

\item[(i)] If $k_1=(p^{4s}-4p^{2s}-5)/2$, we use formula (\ref{main_formula}) from Theorem \ref{Main} and obtain that $8 \sum_{i=1}^8 B_i^2-(\sum_{i=1}^8  B_i)^2=0$. Hence we have $$B_1=B_2=\cdots=B_8=\frac{p^{2s}(p^{2s}-1)}{2}.$$ According to Theorem \ref{Ma2}, since $B_i=\frac{p^{2s}(p^{2s}-1)}{2} \equiv 0 \pmod{p-1}$, we have $h_i \notin D$, $i=1, 2 \cdots, 8$, that is, $D\cap H=\emptyset$.  On the other hand, by applying Proposition \ref{Ma1} to the subgroup $H$, we see that $|D\cap H|=2$ or $8$ (here $\pi=3$, $\theta=-1$, and  $\beta_1=1$), a contradiction.

\medskip

\item[(ii)] If $k_1=(p^{4s}+4p^{2s}-5)/2$, we have $8 \sum_{i=1}^8 B_i^2-(\sum_{i=1}^8  B_i)^2=-48p^{4s}<0$, thus this case will not occur.

\end{itemize}

 Hence  there does not exist a  $(9p^{4s},  r(3p^{2s}+1),-3p^{2s}+r^2+3r,r^2+r)$-PDS with $r=(3p^{2s}-5)/2$. \qed

\begin{theorem}
Let $G$ be an abelian group of order $9p^{4s}$, where $p$ is a prime number  with $p \equiv 1 \pmod 4$, and $s$ is an odd positive integer.  Then there does not exist a  $(9p^{4s},  r(3p^{2s}+1),-3p^{2s}+r^2+3r,r^2+r)$-PDS  $D$ in $G$ when $r=(3p^{2s}-3)/4$ or $r=(3p^{2s}+5)/4$.
\end{theorem}

{\bf Proof:} We use the same notation as in the previous theorem. First assume $r=(3p^{2s}-3)/4$.

  Let $N$ be the subgroup of $G$ of order $p^{4s}$. Applying Proposition \ref{Ma1} to the subgroup $N$,  we have $\pi=p^{2s}$,  $\theta=-1$, $\beta_1=(p^{2s}-3)/2$,  and 
$$k_1= |N \cap D|=\frac{3p^{4s}-3}{4}\; \;\mbox{or} \;\;\frac{p^{4s}+2p^{2s}-3}{4}.$$
Now we discuss the existence of $D$  in two cases based on the possible size of $N \cap D$. 

\begin{itemize}

\item[(i)] If $k_1=(3p^{4s}-3)/4$, using formula (\ref{main_formula}) from Theorem \ref{Main}, we obtain $\sum_{i=1}^8 B_i^2=-3p^{4s}(p^{2s}-1)/2<0$, thus this case will not occur.

\item[(ii)] If $k_1=(p^{4s}+2p^{2s}-3)/4$,  we have  $8 \sum B_i^2-(\sum B_i)^2=0$, and hence 
 $$B_1=B_2=\cdots=B_8=\frac{p^{2s}(p^{2s}-1)}{4}.$$ When $p \equiv 1 \pmod 4$ and $s$  is an odd positive integer, we have  $$B_i =\frac{p^{2s}(p-1)(p^{2s-1}+p^{2s-2}+\cdots+p+1)}{4}\equiv \frac{(p-1)(2s)}{4}\equiv \frac{p-1}{2} \not\equiv 0,1 \pmod {p-1}.$$ Clearly $B_i \not\equiv 0,1 \pmod{p-1}$ contradicts the result of Theorem \ref{Ma2}.

\end{itemize}

Hence  there does not exist a  $(9p^{4s},  r(3p^{2s}+1),-3p^{2s}+r^2+3r, r^2+r)$-PDS with $r=(3p^{2s}-3)/4$. 

\bigskip
 Next we assume that $r=(3p^{2s}+5)/4$. As in the previous case, we apply Proposition \ref{Ma1} to the subgroup $N$. Then we get $\pi=p^{2s}$, $\theta=-1$, $\beta_1=(p^{2s}+5)/2$, and 
 $$k_1=\frac{3p^{4s}+8p^{2s}+5}{4}\; \;\mbox{or} \;\;\frac{p^{4s}-6p^{2s}+5}{4}.$$

\begin{itemize}

\item[(i)] If $k_1=(3p^{4s}+8p^{2s}+5)/4$, using formula (\ref{main_formula}) from Theorem \ref{Main},  we have 
$\sum_{i=1}^8  B_i^2=-p^{4s}(3p^{2s}+5)/2<0$, thus this case will not occur.

\item[(ii)] If $k_1=(p^{4s}-6p^{2s}+5)/4$, we have $8 \sum B_i^2-(\sum B_i)^2=0$, 
and hence $$B_1=B_2=\cdots=B_8=\frac{p^{2s}(p^{2s}+3)}{4}.$$

 When $p \equiv 1 \pmod 4$ and $s$  is an odd positive integer, we have  
 $$B_i= \frac{p^{2s}(p-1)(p^{2s}+\cdots+p+1)}{4}+p^{2s}\equiv \frac{p-1}{2}+1 \not\equiv 0, 1 \pmod {p-1}.$$Clearly $B_i \not\equiv 0,1 \pmod{p-1}$ contradicts the result of Theorem \ref{Ma2}. 
\medskip

\end{itemize}

Hence there does not exist a  $(9p^{4s},  r(3p^{2s}+1),-3p^{2s}+r^2+3r, r^2+r)$-PDS with $r=(3p^{2s}+5)/4$.

\qed

  Now we summarize the results of this section in the following theorem:

\begin{theorem}\label{summary}

\begin{itemize}

\item[{\rm (i)}]
 Let $G$ be an abelian group of order $9m^4$, where $m$ is a positive odd integer and $\gcd(3,m)=1$.  If there exist a  $(9m^4,  r(3m^2+1),-3m^2+r^2+3r,r^2+r)$ negative Latin square type PDS in $G$ with $r\le \frac{3m^2-1}{2}$,  then $$r=\frac{3m^2-3}{4},\; \frac{3m^2+1}{4},\; \frac{3m^2+5}{4}, \; \frac{3m^2-5}{2}, \; \frac{3m^2-3}{2}, \; \mbox{ or }\; \frac{3m^2-1}{2}.$$ 

\item[{\rm (ii)}] Furthermore, if $m=p^s$, where $p$ is a prime number with $p \equiv 1 \pmod 4$, and $s$ is an odd positive integer, under the same assumption as in {\rm (i)}, we have $$r=\frac{3p^{2s}+1}{4},\; \frac{3p^{2s}-1}{2}, \; \mbox{ or }\; \frac{3p^{2s}-3}{2}.$$
\end{itemize}
\end{theorem}

\bigskip

Here we want to mention that there are known constructions of negative Latin square type PDSs with some of the $r$ values from Theorem \ref{summary}:

In 2010, Polhill \cite{Polhill2010} constructed Paley type partial difference sets (negative Latin square type PDS with $r=\frac{3m^2-1}{2}$)  in certain abelian
groups of order  $9m^4$, where $m$ is a positive odd number $>1$.  In 2009, Polhill \cite{Polhill2009} (Theorems 3.1 and 3.2) constructed negative Latin square type PDSs in $G'=\mathbb{Z}_3^2\times \mathbb{Z}_p^{4s}$, where $p$ is an odd prime, with parameters $(9p^{4s}, r(3p^{2s}+1), -3p^{2s}+r^2+3r, r^2+r)$, $r=(3p^{2s}-1)/2$ and $r=(3p^{2s}-3)/2$.  Theorem \ref{summary} shows that Polhill's constructions cover most possible negative Latin square type PDS parameter sets in abelian groups of order $9p^{4s}$ when $p$ is  a prime number with $p \equiv 1 \pmod 4$, and $s$ is an odd positive integer.
\bigskip

{\bf Open problem:} Does there exist a $(9p^{4s}, r(3p^{2s}+1), -3p^{2s}+r^2+3r, r^2+r)$-PDS with $r=(3p^{2s}+1)/4$ when $p$ is a prime number with $p \equiv 1 \pmod 4$, and $s$ is an odd positive integer? We attempted to use the structural information one can derive from the proof of Theorem \ref{Main} to construct such PDS but did not succeed.

\end{document}